\documentclass[12pt,reqno]{amsart}

\usepackage{amsbsy}

\usepackage[colorlinks]{hyperref}
\hypersetup{
linkcolor=blue,          % color of internal links (change box color with linkbordercolor)
citecolor=green,        % color of links to bibliography
}
\usepackage[nobysame,abbrev,alphabetic]{amsrefs}

\usepackage{enumerate}
\usepackage{amssymb}
\usepackage{amsmath}
\usepackage{amscd}
\usepackage{amsthm}
\usepackage{amsfonts}
\usepackage{hyperref}
\usepackage{mathrsfs}

\newtheorem{thm}{Theorem}
\newtheorem*{thm*}{Theorem}
\newtheorem*{remark*}{Remark}

\theoremstyle{definition}

\theoremstyle{definition}

\theoremstyle{definition}
\theoremstyle{definition}
\theoremstyle{definition}

\newcommand{\sm}{\smallsetminus}

\newcommand{\R}{{\mathbb{R}}}

\newcommand{\Z}{{\mathbb{Z}}}

\newcommand {\ignore}[1]  {}

\newcommand{\SL}{\operatorname{SL}}

\newcommand{\onto}{\xymatrix{\ar@{>>}[r]&}}
\newcommand{\da}[4]{\xymatrix{#1 \ar@<.5ex>[r]^{#2} \ar@<-.5ex>[r]_{#3} & #4}}
\newif\ifdraft\drafttrue
\draftfalse

\begin{document}
\title{Counterexamples to 
a conjecture of Woods}
\author{Oded Regev}
\address{Computer Science Department, Courant Institute of Mathematical Sciences, New York
 University.}
\author{Uri Shapira}
\address{Dept. of Mathematics, Technion, Haifa, Israel
%{\tt ushapira@tx.technion.ac.il}
}
\author{Barak Weiss}
\address{Dept. of Mathematics, Tel Aviv University, Tel Aviv, Israel
%{\tt barakw@post.tau.ac.il}
}

\begin{abstract}
A conjecture of Woods from 1972 is disproved. 
\end{abstract}

\maketitle

A lattice in $\R^d$ is called \emph{well-rounded} if its shortest
nonzero vectors span $\R^d$, is called \emph{unimodular} if its
covolume is equal to one, and the \emph{covering radius} of a lattice
$\Lambda$ is the least $r$ such that $\R^d = \Lambda+B_r$, where $B_r$
is the closed Euclidean ball of radius $r$. Let $N_d$ denote the
greatest value of the covering radius over all well-rounded unimodular
lattices in $\R^d$. In \cite{Woods_1972}, A.\ C.\ Woods conjectured that
$N_d = \sqrt{d}/2$, i.e., that the lattice 
$\Z^d$ realizes the largest covering radius among well-rounded unimodular
lattices. Moreover, Woods proved this statement for $d \leq 6$. In
\cite{McMullenMinkowski}, McMullen proved that Woods's
conjecture implies a celebrated conjecture of Minkowski. Spurred by
this result, Woods's conjecture has 
been proved for $d \leq 9$ by Hans-Gill, Kathuria, Raka, and Sehmi
(see~\cite{Chandigarh9} and references therein), thus yielding
Minkowski's conjecture in those dimensions. 
In this note we prove:
\begin{thm*}\label{thm: main} There is $c>0$ such that $N_d >c
  \frac{d}{\sqrt{\log d}}$. 
For all $d \geq 30,$ $N_d > \frac{\sqrt{d}}{2}$.
\end{thm*}

\begin{proof}
%Thus $C(\Z^m)=m$ and Woods's conjecture is that for all
%well-rounded 
%unimodular lattices, $C(\Lambda) \leq \dim \Lambda$. 
Our examples will all be of the form
%\eq{eq: form}{
$$
\Lambda=\alpha_1 
\Lambda_1 \oplus \alpha_2 \Z^m
$$
%}
for some choices of $\Lambda_1,
\alpha_1, \alpha_2, m$. It will be more convenient to work with the
quantity $C(\Lambda) = 
4r(\Lambda)^2$, where $r(\Lambda)$ is the covering radius of
$\Lambda$. 
Clearly $C(\alpha \Lambda) = \alpha^2 C(\Lambda)$, and the Pythagorean
theorem shows that $C(\Lambda_1 \oplus \Lambda_2) =  
C(\Lambda_1)+C(\Lambda_2)$. Let $\lambda_1(L)$ denote the
length of the shortest nonzero vector of $L$, and suppose
$\Lambda_1, \Lambda_2$ are 
well-rounded. If the $\alpha_i$ satisfy $\lambda_1(\alpha_1 \Lambda_1) =
\lambda_1(\alpha_2\Lambda_2)$, then $\alpha_1 \Lambda_1 \oplus
\alpha_2 \Lambda_2$ % as  in \equ{eq: form} 
is well-rounded. Moreover, there is a unique choice of 
$\alpha_i$ for  which it is also unimodular. Namely, if
$\Lambda_1$ is well-rounded and unimodular of dimension $n$, and 
$\Lambda_2 = \Z^m$, in order for
$\Lambda$ to be well-rounded and unimodular we must take 
$\alpha_1 = \lambda^{-\frac{m}{n+m}}$ and $\alpha_2 =
\lambda^{\frac{n}{n+m}}$, where $\lambda = \lambda_1(\Lambda_1)$. Thus 
$$%\eq{eq: formula}{
C(\Lambda) = C(\Lambda_1)  \, \lambda^{-\frac{2m}{n+m}}  +
  m \, \lambda^{\frac{2n}{n+m}}. 
$$%}

For each $d > 3$, let $m =\left \lfloor \frac{d}{\log d} \right \rfloor, \, n =
d-m$. Let $\Lambda_1$ be any lattice in $\R^{n}$ for which
$\lambda_1$ is maximal, that is, $\Lambda_1$ is a lattice giving the
densest lattice packing in dimension $n$. Although $\Lambda_1$ is
only known in very few dimensions, it is a well-known result of
Minkowski (see \cite{GL}*{Chapter 2} or \cite{ConwaySloane}*{\S 1.1.5})
that there is $c_1>0$ such that for all $n$, 
$$\lambda = \lambda_1(\Lambda_1) \geq c_1 \sqrt{n}.$$ 
Recall that a
lattice $L_0$ is called \emph{critical} if the function
$L \mapsto \lambda_1(L)$, considered as a function on the
space of unimodular lattices, attains a local maximum at
$L_0$. Then clearly $\Lambda_1$ is critical, and a theorem of
Voronoi (whose proof is not difficult; see, e.g.,~\cite{GL}*{Chapter
  6}) implies that $\Lambda_1$ is 
well-rounded. Now let $\alpha_1, \alpha_2$ be the unique positive
numbers for which $\Lambda = \alpha_1 \Lambda_1 \oplus \alpha_2
\Z^{m}$ is 
well-rounded and unimodular.
Then
%$$
%$
%\alpha_2 = \lambda^{\frac{n}{d}} \geq c_2 \, \sqrt{d}
%$$
%$
%for some $c_2>0$ and hence 
$$C (\Lambda) \geq m\, \lambda^{\frac{2n}{m+n}} \geq 
c_2 \, m \, n^{\frac{n}{d}}  \geq 
%c_3 \, \frac{d}{\log d} \left(\frac{\log d
%  -2}{\log d}\right)^{\frac{\log d -2}{\log d}} d^{1-\frac{2}{\log d}}
%\geq  
c_3 \frac{d^2}{\log d}$$
%the covering radius of $\alpha_2 \Z^{d_2}$ is at least
%$\alpha_2 \frac{\sqrt{d_2}}{2} \geq 
%c_3 \, d^{\frac{1}{2}-
 % \frac{1}{2\log d}} \, d_2^{1/2}
%\geq c_4 \, \frac{d}{\log d}$ 
for positive 
$c_2, c_3$, and the first assertion follows. 

Taking $\Lambda_1$ to be the laminated lattice $\Lambda_{15}$ (see
\cite{ConwaySloane}*{Chapter 6}), we have\footnote{Actually $C(\Lambda_{15})=7 \cdot 2^{\frac{2}{5}}$, as was shown in \cite{Mathieu_computations}; we only need the lower bound.}
$C(\Lambda_1) \geq 
%\frac{14}{512^{1/15}}, 
7 \cdot 2^{\frac{2}{5}}, \, 
\lambda = 
2^{\frac{7}{10}}, \,
%\frac{2}{512^{1/30}}, 
n=15$ and so 
$$
C(\Lambda) \ge 
7 \cdot 2^{\frac{2}{5}} \,
%\frac{14}{512^{1/15}} \cdot
\left(%\frac{2}{512^{1/30}} 
2^{\frac{7}{10}} \right)^{-\frac{2m}{15+m}} + m \, \left(
    %\frac{2}{512^{1/30}} 
2^{\frac{7}{10}}
\right)^{\frac{30}{15+m}},
$$
which is greater than $d=m+15$ for all $m \geq 15$, as can be seen by
evaluating the function at $m=15$ and showing that its derivative is
positive for $m \geq 15$. Note that
$\Lambda_{15}$ is generated by its shortest nonzero vectors and is in
particular well-rounded. See \cite{ConwaySloane}*{Chapter 6} or \cite{Barnes}. 
%Barnes-Wall lattice in dimension 16 we
%have $C(\Lambda_1)= 6\sqrt{2}, \lambda = 2^{3/4}, d =16$ (see
%\cite[Chapter 4.10]{ConwaySloane}) and so 
%$$
%C(\Lambda) = 6\sqrt{2} \cdot 2^{\frac{-6m}{64+4m}} + m \cdot 2^{\frac{96}{64+4m}},
%$$
%which is greater than $m+16$ for all $m \geq 17.$ This proves
%the first assertion 
\end{proof}

\begin{remark*}
A similar construction using different choices of $\Lambda_1$ (instead
of $\Lambda_{15}$) also
works, giving slightly weaker results: the Leech lattice (for $d \geq 38$), the `shorter Leech lattice' $O_{23}$ (for $d
\geq 36$), the 16-dimensional Barnes-Wall lattice (for $d \geq 33$),
and the laminated lattice $\Lambda_{23}$ 
will work (for $d \geq 31$). We are grateful to M. Dutour-Sikiri\'c for
suggesting the use of various lattices, and in particular the
laminated lattice $\Lambda_{15}$ for this 
problem. The covering radii of these and other lattices can be
computed using his publicly available program \cite{Polyhedral}. 
%
%Taking $\Lambda_1$ to be the laminated lattice $\Lambda_{23}$ in
%dimension 23, we have
%$C(\Lambda_1 )= 15 \cdot 2^{2/23}$ and $\lambda = 2^{22/23}$, and so 
%$C(\Lambda) =  15 \cdot 2^{2/23} 2^{-\frac{44x}{23(23+x)}} + x \,
%2^{\frac{22}{23+x}}$ which is greater than $23+x$ when $x \geq  8$. 
\ignore{
Taking $\Lambda_1$ to be the Leech lattice we have $c = \sqrt{2}, \lambda
=2, d =24$ (see \cite[Chapter 4.10]{ConwaySloane}) and so 
$$
C(\Lambda) = 2^{\frac{48}{24+m}}(m+2),
$$
which is greater than $\frac{\sqrt{m+24}}{2}$ for all $m \geq 14.$ This proves
the first assertion for $d \geq 38$. A
a similar computation with $\Lambda_1$ equal to the Barnes-Wall
lattice in dimension 16 
gives a similar example in dimension $d \geq 33$. 
}
\end{remark*}

{\bf Acknowledgements:}
We are grateful to Mathieu Dutour-Sikiri\'c and Curt McMullen for useful
comments and suggestions. OR was supported by the Simons Collaboration
on Algorithms and Geometry and by the National Science Foundation
(NSF) under Grant No.~CCF-1320188. US was supported by ISF grant
357/13. BW was supported by ERC starter grant
DLGAPS 279893. 
%Any opinions, findings, and conclusions or
%recommendations expressed in this material are those of the authors
%and do not necessarily reflect the views of the NSF. 
%\end{document}
%\end{appendix}
\bibliographystyle{alpha}
\bibliography{elonbib}
\end{document}